\documentclass[12pt]{article}
\usepackage{amssymb,latexsym,amsmath,amsthm,amsfonts,txfonts,rotating}

\def\im{{\rm im}}
\newcommand{\nc}{\newcommand}
\renewcommand{\H}{{\mathbb H}}
\nc{\Hom}{\mathop{\rm Hom}\nolimits}
\nc\vv{\mathbf v}
\nc\w{\mathbf w}
\nc\calV{I}
\nc\calE{E}
\nc\gG{{\mathcal G}}
\nc\g{\mathfrak g}
\nc\gl{{{\mathfrak g}{\mathfrak l}}}
\nc{\GL}{{\rm GL}}
\nc{\SL}{{\rm SL}}
\nc{\PGL}{{\rm PGL}}
\nc{\G}{{\rm G}}
\nc{\T}{{\mathbb T}}
\nc{\F}{{\mathbb F}}
\nc{\mF}{{\mathcal F}}
\nc{\mC}{{\mathcal C}}
\nc{\calM}{{\cal M}}
\nc{\calP}{{\cal P}}
\nc{\calT}{{\cal T}}
\nc{\N}{{\mathbb N}}
\nc{\Q}{{\mathbb Q}}
\nc{\J}{{\cal J}}
\nc{\higgsinfty}{${{\rm Higgs}_\infty}\ $}
\nc{\lrow}{\longrightarrow}
\nc{\Dol}{{\op{Dol}}}
\nc{\DR}{{\op{DR}}}
\nc{\B}{{\op{B}}}
\nc{\op}[1]{\mathop{\mathchoice{\mbox{\rm #1}}{\mbox{\rm #1}}
{\mbox{\rm \scriptsize #1}}{\mbox{\rm \tiny #1}}}\nolimits}
\nc{\down}[1]{{\phantom{\scriptstyle #1}
\hbox{$\left\downarrow\vbox to
    9.5pt{}\right.\nulldelimiterspace=0pt \mathsurround=0pt$}
\raisebox{.4ex}{$\scriptstyle #1$}}}
\nc{\Qu}{{\mathbb Q}}
\nc{\calC}{{\cal C}}
\nc{\calD}{{\cal D}}
\nc{\stack}[2]{
{\begin{array}{c}
\scriptstyle #1 \\ \scriptstyle #2 \end{array}} }

\nc{\calB}{{\cal B}}
\nc{\C}{{\mathbb C}}
\nc{\R}{{\mathbb R}}
\nc{\calK}{{\cal K}}
\nc{\calR}{{\cal R}}
\nc{\Proj}{{\mathbb P}}
\nc{\Est}{E_{\op{st}}}
\nc{\Hst}{H_{\op{st}}}
\nc{\Z}{{\mathbb Z}}
\nc{\Nat}{{\mathbb N}}
\nc{\E}{{\mathbb E}}
\nc{\U}{{\rm U}}
\nc{\SU}{{\rm SU}}
\nc{\dbar}{{\overline{\partial}}}
\nc{\K}{\mathbb K}
\nc{\Dirac}{{D\!\!\!\! \slash}} 
\nc{\Ha}{{\cal H}}
\nc{\Hy}{{\mathbb H}}
\nc{\V}{{\mathbb V}}
\nc{\M}{{\mathcal M}}
\nc{\I}{{\mathbb I}}
\renewcommand{\J}{{\mathbb J}}
\nc{\A}{{\mathcal A}}
\nc{\uPhi}{{\mathbf \Phi}}
\nc{\cM}{{{\overline\M}}}
\nc{\bM}{{\mathbb M}}
\nc{\calO}{{\cal O}}
\nc{\beq}{\begin{eqnarray}}
\nc{\eeq}{\end{eqnarray}}
\nc{\bes}{\begin{eqnarray*}}
\nc{\ees}{\end{eqnarray*}}
\nc{\oper}[1]{\mathop{\mathchoice{\mbox{\rm #1}}{\mbox{\rm #1}}
{\mbox{\rm \scriptsize #1}}{\mbox{\rm \tiny #1}}}\nolimits}
\nc{\operlimits}[1]{\mathop{\mathchoice{\mbox{\rm #1}}{\mbox{\rm #1}}
{\mbox{\rm \scriptsize #1}}{\mbox{\rm \tiny #1}}}}
\nc{\Coeff}{\operlimits{Coeff}}
\nc{\Left}[1]{\hbox{$\left#1\vbox to
    10.5pt{}\right.\nulldelimiterspace=0pt \mathsurround=0pt$}}
\nc{\Right}[1]{\hbox{$\left.\vbox to
    10.5pt{}\right#1\nulldelimiterspace=0pt \mathsurround=0pt$}}
\nc{\LEFT}[1]{\hbox{$\left#1\vbox to
    15.5pt{}\right.\nulldelimiterspace=0pt \mathsurround=0pt$}}
\nc{\RIGHT}[1]{\hbox{$\left.\vbox to
    15.5pt{}\right#1\nulldelimiterspace=0pt \mathsurround=0pt$}}

\newtheorem{theorem}{Theorem}[section]

\newtheorem{corollary}[theorem]{Corollary}
\newtheorem{conjecture}[theorem]{Conjecture}

\newtheorem{problem}[theorem]{Problem}
\theoremstyle{remark}
\newtheorem{remark}[theorem]{Remark}
\newtheorem{example}[theorem]{Example}

\theoremstyle{definition}

 \numberwithin{equation}{section}

\begin{document} 

\title{S-duality in Hyperk\"ahler Hodge Theory}

\author{ Tam\'as Hausel\footnote{ University of Oxford \& University of Texas at Austin, \tt hausel@maths.ox.ac.uk}}
\date{}
  \maketitle
{\begin{center}{\it To Nigel Hitchin for his 60th birthday.}\end{center}}

\begin{abstract} Here we survey questions and results on  the Hodge theory of hyperk\"ahler quotients, motivated by certain S-duality considerations in string theory. The problems include $L^2$ harmonic forms, Betti numbers and mixed Hodge structures on 
the moduli spaces of Yang-Mills instantons on ALE gravitational instantons, magnetic monopoles on $\R^3$ and Higgs bundles
on a Riemann surface. Several of these spaces and their hyperk\"ahler metrics
were constructed by Nigel Hitchin and his collaborators. 
\end{abstract}

\section{Introduction}

In this paper we survey the motivations, related results and progress made towards the following problem, raised by Hitchin in 1995:

\begin{problem}\label{problem} What is the space of $L^2$ harmonic forms on the moduli space
of Higgs bundles on a Riemann surface?
\end{problem}

The moduli space $\M^d_{\Dol}(\SL_n)$ of stable rank $n$ Higgs bundles with fixed determinant of degree $d$  on a Riemann surface was introduced and studied in \cite{hitchin-self}, \cite{simpson-nonabelian} and \cite{nitsure}. The Betti numbers  of this space for $n=2$ were determined in \cite{hitchin-self} while for $n=3$ in \cite{gothen}.  The above problem raised two new directions to study.
First is the Riemannian geometry of $\M^d_{\Dol}(\SL_n)$, or more precisely the asymptotics of the natural hyperk\"ahler metric, and its connection with Hodge theory. The second one, which can  be considered the  topological side of Problem~\ref{problem}, is to determine the  intersection form on the middle dimensional compactly supported cohomology of $\M^d_{\Dol}(\SL_n)$. While the first question seems still out of reach, although we will report on some modest progress below, the second is more approachable and we offer a conjecture at the end of this survey. 

Problem~\ref{problem} was motivated by $S$-duality conjectures emerging 
from the string theory literature about Hodge theory on certain
hyperk\"ahler moduli spaces, which are close relatives of $\M^d_{\Dol}(\SL_n)$.  

In the physics literature $S$-duality stands for a {\em strong-weak duality} between two quantum field theories. The interest from the physics point of view is that it gives a tool to study physical theories with a large coupling constant via a conjectured equivalence with a theory  with a small coupling constant where perturbative methods give a good understanding. The $S$-duality conjecture relevant for us is based on the Montonen-Olive electro-magnetic duality proposal from 1977 in four dimensional Yang-Mills theory \cite{montonen-olive}. It was noted 
in \cite{olive-witten} that this duality proposal is more  likely to hold in a  supersymmetric version of the theory, and in \cite{osborn} it was argued that  $N=4$ supersymmetry is a good candidate.  Hyperk\"ahler Hodge theory is relevant in $N=4$ supersymmetry as the space of differential forms on a hyperk\"ahler manifold possesses an action of the $N=4$ supersymmetry algebra via the various operators 
in hyperk\"ahler  Hodge theory.

In this paper our interest lies in the mathematical predictions of such
$S$-duality conjectures in physics.
In 1994 Sen \cite{sen}, using   $S$-duality arguments in $N=4$ supersymmetric Yang-Mills theory, predicted the dimension of the spaces 
$\Ha^d(\widetilde{M}_k^0)$ of $L^2$ harmonic $d$-forms on the universal cover $\widetilde{M}_k^0$  of the hyperk\"ahler moduli space $M_k^0$ of certain ${\rm SU}(2)$ magnetic monopoles on $\R^3$. In the interpretation of \cite{sen} the $L^2$ harmonic forms on 
$\widetilde{M}_k^0$ can be thought of as bound states of the theory, and the conjectured $S$-duality implies an action of ${\rm SL}(2,\Z)$ on $\bigoplus_k \Ha^*(\widetilde{M}_k^0)$. 
By further physical arguments Sen managed to predict this representation of ${\rm SL}(2,\Z)$ completely, implying the following

\begin{conjecture}  \label{sen} The dimension of the space of $L^2$ harmonic forms on $\widetilde{M}_k^0$ is 
$$\dim\left(\Ha^d(\widetilde{M}_k^0)\right)=\left\{\begin{array}{ll} 0 & d\neq mid \\ \phi(k) & d=mid  \end{array} \right.,$$ where $\phi(k)= \sum_{i=1}^{k} \delta_{1 (i,k)}$ is the Euler $\phi$ function, and $mid=2k-2$ is half of the  dimension of $\widetilde{M}^0_k$.
 \end{conjecture}

Similar $S$-duality arguments led Vafa and Witten \cite{vafa-witten}
to get a conjecture on the space of $L^2$ harmonic forms on a certain smooth completion ${M}^{k,c_1}_\phi $, constructed in \cite{kronheimer-nakajima, nakajima-quiver2}, of the moduli space of ${\rm U}(n)$ Yang-Mills instantons 
of first Chern class $c_1$, energy $k$ and framing $\phi$ on one of Kronheimer's  ALE 
spaces, which are $4$-dimensional  complete hyperk\"ahler manifolds, 
with an asymptotically locally Euclidean metric.  

\begin{conjecture}\label{vafawitten}   The dimension of the space of $L^2$ harmonic forms on $M^{k,c_1}_\phi$ is $$\dim\left(\Ha^{d}(M^{k,c_1}_\phi)\right)=\left\{ \begin{array}{ll} 0 & d\neq mid\\
\dim\left(\im(H_{cpt}^{mid}(M^{k,c_1}_\phi)\rightarrow H^{mid}(M^{k,c_1}_\phi))\right) & d=mid\end{array} \right. ,$$ where $mid$ now denotes half of the dimension of $M^{k,c_1}_\phi$. 
\end{conjecture}
\noindent The paper \cite{vafa-witten} further argues that Conjecture~\ref{vafawitten} 
implies, via the work of Nakajima \cite{nakajima-quiver2} and 
Kac \cite{kac-book},
that \beq \label{modular} Z_{\phi}(q)=\sum_{c_1,k} q^{k-c/24} \dim\left(\Ha^{mid}(M^{k,c_1}_\phi)\right)\eeq is a modular form, which, as was speculated in \cite{vafa-witten}, might be a consequence of  
$S$-duality. 

This paper will introduce the reader to various mathematical aspects of  these three problems and offer
 mathematical techniques and results relating to them. 

\begin{paragraph} {\bf Acknowledgment.}  This paper is a write-up of 
the author's talk at the Geometry Conference in Honour of Nigel Hitchin in Madrid in September 2006. Problem~\ref{problem} was raised by Nigel Hitchin in 1995, then the author's PhD supervisor, as a project for the author's PhD thesis.  This survey paper would like to show the impact of this modestly looking question on the author's subsequent research.   The author's research has been supported by a Royal Society University Research Fellowship, NSF grant DMS-0604775 and an Alfred Sloan Fellowship. The visit to Madrid was supported by a Royal Society International Joint Project between the UK and Spain. 
\end{paragraph}

\section{Hyperk\"ahler quotients}
A Riemannian manifold $(M,g)$ is hyperk\"ahler if it is K\"ahler with 
respect to three integrable complex structures $I,J,K\in\Gamma({\rm End}(TM))$, which satisfy $I^2=J^2=K^2=IJK=-1$, with K\"ahler forms $\omega_I$, $\omega_J$ and $\omega_K$. Known  compact examples are scarce, see e.g. \cite[\S7]{joyce-book}. Non-compact complete examples however are much more abundant. This is mostly because there is a widely applicable\footnote{Some colleagues even suggest, due to the success of this construction, that HyperK\"ahLeR is in fact just a pronouncable version of the acronym HKLR.} {\em hyperk\"ahler quotient construction}, due to Hitchin-Karlhede-Lindstr\"om-Ro{\v{c}}ek 
\cite{hitchin-etal}. 
The construction itself is an elegant quaternionization of the Marsden-Weinstein symplectic (or more precisely K\"ahler) quotient construction (see \cite[Chapter 8]{mumford-etal} for an introduction for the latter). 

Let  ${\mathbb M}$ be a
hyperk\"ahler manifold,   $\gG$ a 
Lie group, with Lie algebra $\g$,  and assume 
$\gG$ acts on  ${\mathbb M} $ preserving the hyperk\"ahler structure (i.e. it acts by triholomorphic isometries). 
Let us  further assume that we have  moment maps $\mu_I:\bM\to \g^*$, $\mu_J:\bM\to \g^*$ and $\mu_J:\bM\to \g^*$ with respect to
the symplectic forms $\omega_I$, $\omega_J$ and $\omega_K$ respectively. 
We combine them into a single hyperk\"ahler moment map: 
$$\mu_\H=(\mu_I,\mu_J,\mu_K): \bM \to \R^3\otimes \g^*.$$
One takes
$\xi\in \R^3\otimes (\g^*)^\G$ and constructs
the {\em hyperk\"ahler quotient} at level $\xi$ by:
$$\bM /\!/\!/\!/_\xi \gG:=\mu_\H^{-1}(\xi)/\gG.$$
The main result of \cite{hitchin-etal} is that the natural Riemannian metric
on the smooth points of this quotient is hyperk\"ahler.

Now we list three important examples of this construction, where the
original hyperk\"ahler manifold $\bM$ and Lie group $\gG$ are both  infinite dimensional.

\subsection{Moduli of Yang-Mills instantons on $\R^4$}
\label{YM}

Here we follow \cite[I Example 3.6]{hitchin-book},  compare also with \cite{atiyah-YMbook}.

Let  $\G$ be a compact connected Lie group, which will be $\U(n)$ or $\SU(n)$ in this paper. Let $P\to \R^4$ be a $\G$-principal bundle over $\R^4$. 
Let  $\bM$ be the space of $\G$-connections  $A$  on P of class  $C^\infty$, such that the energy  $$\left|\int_{\R^4} tr(F_A\wedge *F_A) \right| <\infty
 $$ is finite. 
Write $$A=A_1dx_1+A_2dx_2+A_3dx_3+A_4dx_4$$ in a fixed gauge,
where $A_i\in \Omega^0(\R^4, \rm{ad} (P))$. 
 Let $\gG=\Omega(\R^4,Ad(P))$  be the gauge group of $P$. An element $g\in \gG$ acts on $A\in \bM$  by the formula $g(A)=g^{-1}Ag+g^{-1}dg $, preserving the hyperk\"ahler structure. 
  One finds that the hyperk\"ahler moment map equation $$\mu_\H(A)=0 \Leftrightarrow F_A=*F_A$$ is just the
self-dual Yang-Mills equation. Define the hyperk\"ahler quotient
  $\M(\R^4,P)=\mu_\H^{-1}(0)/\gG$, the moduli
space of finite energy self-dual Yang-Mills instantons on $P$. By its construction it has a natural
hyperk\"ahler metric. 

 Similar construction \cite{kronheimer-nakajima} for $\G=\U(n)$  yields a hyperk\"ahler metric on moduli spaces of $\U(n)$ Yang-Mills instantons on certain four dimensional complete hyperk\"ahler manifolds, the ALE spaces of Kronheimer \cite{kronheimer}.  These moduli spaces will have natural completions  and various components of them will be the spaces $M^{k,c_1}_\phi$ which were mentioned in the introduction. They  will resurface later
 as examples for Nakajima quiver varieties.

\subsection{Moduli space of magnetic monopoles on $\R^3$}\label{hitchin}
The following construction can be considered as a dimensional reduction of the previous example.  Here we follow \cite[I Example 3.5]{hitchin-book} and \cite{atiyah-hitchin-book}. 
 
 Assume that $\G=\SU(2)$ and the matrices $A_i$ are independent of $x_4$. Then we have 
 $$A\!=\!A_1dx_1\!+\!A_2dx_2\!+\!A_3dx_3$$  a connection on $\R^3$ and  $A_4=\phi\in \Omega^0(\R^3,{\rm ad} P)$ becomes the
  {\em Higgs field}. The gauge group now will be
 $\gG=\Omega(\R^3,{\rm Ad} P)$ and $\bM=\{(A,\phi)\mbox{ +  certain boundary condition}\}$. (The boundary condition is chosen to ensure finite energy.) The gauge group $\gG$ acts on $\bM$ by gauge transformations, preserving the natural hyperk\"ahler metric on $\bM$.
 The corresponding hyperk\"ahler moment map equation
 $$\mu_\H(A,\phi)=0\Leftrightarrow F_A=*d_A\phi$$ is equivalent with the Bogomolny equation.
 
 Now by construction $M=\mu_\H^{-1}(0)/\gG$, the moduli
space of magnetic monopoles on $\R^3$, has a natural
hyperk\"ahler metric. It has infinitely many components $M=\cup_{k=1}^\infty M_k$ labeled by the  magnetic charge $k$ of the monopole. 

$M_k$ is acted upon by $\R^3$ by translations and by $\U(1)$ by rotating the phase of the monopole. The quotient $M^0_k$ is still a smooth complete hyperk\"ahler manifold of dimension $4k-4$, with fundamental group $\Z_k$. We will denote by $\widetilde{M_k^0}$ its universal cover. 
In \cite{atiyah-hitchin-1985} Atiyah and Hitchin  find the hyperk\"ahler metric explicitly  on  the  $4$-manifold $M^0_2$ and subsequently  describe the scattering of two monopoles. 

\subsection{Hitchin moduli space} This example can
be considered as a two-dimensional reduction of \S\ref{YM}. We follow \cite[Section 1]{hitchin-self} and \cite[I Example 3.3]{hitchin-book} .

 Now we assume that $G=\U(n)$ and the matrices $A_i$ in \S\ref{YM} are independent of $x_3$, $x_4$. We have now  the connection
  $A\!=\!A_1dx_1\!+\!A_2dx_2$ on the $\U(n)$ principal bundle $P$ on $\R^2$. We introduce
   $\Phi=(A_3-A_4i)dz\in  \Omega^{1,0}(\R^2,{\rm ad} P\otimes \C)$ the {\em complex Higgs field}. The gauge group now is
  $\gG=\Omega(\R^2,{\rm Ad} P)$, which   acts by gauge transformations on the space $\bM=\{(A,\Phi)\}$ preserving 
  the natural hyperk\"ahler metric on $\bM$.
   The moment map equations \begin{eqnarray*}\mu_\H(A,\Phi)=0 \Leftrightarrow
\begin{array}{c} F(A)=-[\Phi,\Phi^*],\\
d_A^{\prime\prime}\Phi=0.\end{array} \label{equation} 
\end{eqnarray*}
are then
equivalent with Hitchin's self-duality equations. There are no solutions of finite energy on $\R^2$, but as the equations are conformally invariant, we can replace $\R^2$  
with a genus $g$ compact Riemann surface $C$ in the above definitions, and define
$\M(C,P)=\mu_\H^{-1}(0)/\gG$, the Hitchin moduli space, which 
has a natural hyperk\"ahler metric by construction. There are different ways to think about this space  with the different complex structures, which will be explained in \S\ref{diffeo}.

\section{Hodge theory}
\subsection{$L^2$ harmonic forms on complete manifolds}

Let  $M$ be a  complete Riemannian manifold of dimension $n$.  We say that a
smooth differential $k$-form
$\alpha \in \Omega^k(M)$ is harmonic if and only if
$d\alpha=d\!*\! \alpha=0$, where $*:\Omega^k(M)\to \Omega^{n-k}(M)$
is the Hodge star operator.  It is  $L^2$ if and only if $$\int_M \alpha\wedge *\alpha<\infty.$$ We denote by $\Ha^*(M)$  the space of $L^2$ harmonic forms. 

A fundamental theorem  of Hodge theory is the 
  Hodge (orthogonal) decomposition theorem 
\cite[\S32 Theorem 24, \S35 Theorem 26]{derham} : \beq \label{dec}\Omega^*_{L^2}=\overline{d(\Omega^*_{{cpt}})}\oplus \Ha^*
  \oplus 
\overline{\delta(\Omega_{cpt}^*)},\eeq where $\delta$ is the adjoint of 
$d$. When $M$ is compact this implies the celebrated Hodge theorem, which 
says that 
$\Ha^*(M)\cong H^*(M)$ i.e. that there is a unique harmonic representative in every de Rham cohomology class. 
When $M$ is non-compact we only have a topological lower bound.
Namely, the Hodge decomposition theorem  implies that the composite map
  $$H^*_{cpt}(M)\rightarrow \Ha^*(M)\rightarrow H^*(M)$$ is just the forgetful map. (In the compact case these maps are isomorphisms, which gives the Hodge theorem mentioned above.) Thus \beq \label{lower} \im(H^*_{cpt}(M)\to H^*(M))\eeq is a "topological lower bound"
for $\Ha^*(M)$. By Poincar\'e duality the map  $H_{cpt}^*(M)\rightarrow H^*(M)$
is equivalent with the intersection
pairing on $H^*_{cpt}(M)$. 

In the cases most relevant for us $M$ will be a hyperk\"ahler manifold (sometimes orbifold) so $\dim(M)=4k$ and we will additionally have $H^i(M)=0$ for $i>2k$. Therefore
the possible non-trivial image in $\im(H^*_{cpt}(M)\to H^*(M))$ will
be concentrated in the middle $2k$ dimension. (We will use the notation $mid=\dim(M)/2$ for the middle dimension of a manifold.) For such a hyperk\"ahler manifold we denote \beq \label{chitopl2} \chi_{L^2}(M)=\dim\left(\im(H^{mid}_{cpt}(M)\to H^{mid}(M))\right)=\dim\left(\im(H^*_{cpt}(M)\to H^*(M))\right)\eeq
the dimension  of this image. $\chi_{L^2}(M)$ can be thought of either as a "topological lower bound" for $\dim(\mathcal{H}^*(M))$ or the Euler characteristic of topological $L^2$ cohomology. 
\subsection{Results on $L^2$ harmonic forms}
There were few general theorems on describing $\Ha^*(M)$ for a non-compact complete manifold $M$, see however \cite[Introduction]{hausel-hunsicker-mazzeo} for an overview. 
It was thus a surprising development when Sen \cite{sen}, using 
arguments from $S$-duality, managed to predict 
the dimension of $L^2$ harmonic forms on $\widetilde{M}_0^k$ as was
explained in Conjecture~\ref{sen} in the Introduction. 
 In particular, according to Sen's Conjecture~\ref{sen} the space $\Ha^2(\widetilde{M}_2^0)$ should be one dimensional. Using the explicit description \cite{atiyah-hitchin-1985} of the  metric on $\M_2^0$   Sen in \cite{sen} was able to 
  find an explicit $L^2$ harmonic $2$-form, called the {\em Sen $2$-form}, on $\widetilde{M}_2^0$. 
  This was perhaps the strongest mathematical support exhibited 
  for Conjecture~\ref{sen} in \cite{sen}. 
  
  More general mathematical support for Conjecture~\ref{sen} came  in 1996. Segal and Selby in \cite{segal-selby} showed that 
  the intersection form on $H_{cpt}^{mid}(\widetilde{M}_k^0)$ is definite. Moreover they  obtained for the topological  lower bound \eqref{lower} for ${\Ha}^{mid}(\widetilde{M}_k^0)$  $$\chi_{L^2}\left(\widetilde{M}_k^0\right)=\dim\left(H^{mid}(\widetilde{M}_k^0)\right)= \phi(k).$$ This agrees with the predicted dimension of ${\Ha}^{mid}(\widetilde{M}_k^0)$ in Sen's Conjecture~\ref{sen} . 

Motivated by Problem~\ref{problem} and Segal-Selby's topological lower bound for Conjecture~\ref{sen}, \cite{hausel-vanishing} calculated in 1998 that the intersection pairing on the $g$ dimensional space $H_{cpt}^{mid}(\M^1_{\Dol}(\SL_2))$ is trivial, in other words  \beq \label{trivial}\chi_{L^2}\left(\M^1_{\Dol}(\SL_2)\right)=0\eeq for $g>1$.  This thus gave the surprising result that there are no $L^2$ harmonic forms on $\M^1_{\Dol}(\SL_2)$ plainly by topological reasons.  The technique
used in the proof of \eqref{trivial} was imitating Kirwan's proof \cite{kirwan-mumford} of Mumford's conjecture on the cohomology ring of the moduli space of
stable rank $2$ bundles of degree $1$ on the Riemann surface $C$. Therefore the extension of \eqref{trivial} to higher rank Higgs bundle moduli spaces $\M^d_{\Dol}(\SL_n)$  was not straightforward.

  Next advance towards Sen's Conjecture~\ref{sen} came in 2000.  Hitchin in \cite{hitchin-l2} showed  that   $\Ha^d(M)=0$ unless $d=\dim(M)/2$ for a complete 
hyperk\"ahler manifold $M$ of linear growth.  Examples include all our hyperk\"ahler quotients discussed in this paper.  The proofs in \cite{hitchin-l2} use techniques inspired
  by Jost an Zuo's extension \cite{jost-zuo} of ideas of Gromov \cite{gromov}. It is interesting to note that some of the proofs in \cite{hitchin-l2} also exploit the operators in hyperk\"ahler Hodge theory, which are relevant in $N=4$ supersymmetry. Using the symmetries of the Atiyah-Hitchin metric \cite{hitchin-l2} proves Sen's conjecture for $k=2$, i.e. that up to a scalar the only $L^2$ harmonic form on $\widetilde{M}_2^0$ is Sen's $2$-form. 

A more topological approach was introduced in \cite{hausel-hunsicker-mazzeo} in 2004. \cite{hausel-hunsicker-mazzeo}  proves for fibered boundary
 manifolds $M$ \beq \label{middle} \Ha^{mid}(M)\cong\im(IH_{\underline{m}}^{mid}(\overline{M}){\rightarrow} IH_{\bar{m}}^{mid}(\overline{M})),\eeq where $\overline{M}$ is a certain compactification of $M$, dictated by the asymptotics of the fibered boundary metric on $M$. Moreover $IH_{\underline{m}}^{mid}(\overline{M}))$ denotes the intersection cohomology in dimension $mid=\dim(M)/2$ with lower middle perversity $\underline{m}$ and $IH_{\bar{m}}^{mid}(\overline{M}))$ denotes the intersection cohomology in the middle dimension with upper middle perversity $\bar{m}$ of the possibly badly singular (i.e. not necessarily a Witt space) compactification $\overline{M}$. 
To illustrate \eqref{middle} we take  the compactification of $\widetilde{M}_2^0$, which  happens to be the smooth space $\C{\mathbb P}^2$ (with the non-standard orientation),  where the above cohomologies in \eqref{middle} all coincide, giving $\Ha^2(\widetilde{M}_2^0)\cong H^2(\C{\mathbb P}^2)$. This provides a topological explanation for the existence and uniqueness of the Sen $2$-form. 
 
 The assumption that the metric is fibered boundary in \cite{hausel-hunsicker-mazzeo} is fairly restrictive. Among  hyperk\"ahler quotients
 only a few examples satisfy this property (see the discussion in \cite[\S7]{hausel-hunsicker-mazzeo}). Examples include all  ALE gravitational instantons of \cite{kronheimer} and all known ALF (see \cite{cherkis-kapustin-ALF}) and some ALG gravitational instantons (see \cite{cherkis-kapustin-ALG}). 
In general our hyperk\"ahler quotients have some kind of stratified asymptotic behaviour at infinity. For example the metric  on ${M}_k^0$ is fibered boundary only when $k=2$, for higher $k$ it is known to behave differently at different regions of infinity.  The first result, which could handle Hodge theory on Riemannian manifolds with such a stratified behaviour at infinity appeared recently  in a work \cite{carron-qale} by Carron. It proves  for a QALE space $M$ that:
 $$\Ha^{mid}(M)\cong\im(H_{cpt}^{mid}({M}){\rightarrow} H^{mid}({M})).$$
 A QALE space \cite[\S9]{joyce-book} by definition is a certain Calabi-Yau metric on a crepant resolution of $\C^k/\Gamma$, where $\Gamma\subset \rm{SU}(k)$ is a finite subgroup. The asymptotics of the metric on such a QALE space is reminiscent to the asymptotics of the natural hyperk\"ahler metric
 on ${M}^{k,c_1}_\phi$ appearing in the Vafa-Witten Conjecture~\ref{vafawitten}. It is thus reasonable to hope that the Vafa-Witten Conjecture~\ref{vafawitten} will be decided soon. 
 
 As there have been extensive studies starting with \cite{gibbons-manton} and more recently \cite{Bielawski:2007tq}  on the asymptotics of the Riemannian metric on $M_k^0$, it is conceivable that we will have a precise understanding of the asymptotic behaviour of this metric, and in turn the Hodge theory of $L^2$ harmonic forms on $\widetilde{M}^0_k$, perhaps extending techniques from \cite{carron-qale}. Thus one may be optimistic that  Sen's Conjecture~\ref{sen}  will be decided in the foreseeable future. 
 
 Finally, one must admit that the description of the asymptotics of
 the metric at infinity on $\M^d_{Dol}(\SL_n)$ is still lacking, thus calculation of $\Ha^*\left(\M^d_{Dol}(\SL_n)\right)$ is presently hopeless.
 The topological side of Problem~\ref{problem}, that is to determine
 $\chi_{L^2}\left(\M^d_{\Dol}(\SL_n)\right)$, when $(d,n)=1$,  is more reasonable. After introducing a new arithmetic technique to study Hodge structures
 on the cohomology of our hyperk\"ahler manifolds, we will be able to offer a general conjecture on the intersection form on Higgs moduli spaces, in particular that \eqref{trivial} holds for any $n$. 

\section{Mixed Hodge theory}
\label{arithmetic}
As explained above there have been some limited successes of calculating $\Ha^*(M)$ for a hyperk\"ahler quotient and understanding its relation to the cohomology $H^*(M)$ or more generally the cohomology of an appropriate compactification $H^*(\bar{M})$.
 Another extension of Hodge theory yields some different and in some ways 
 more detailed insight into the cohomology of our 
 hyperk\"ahler quotients. This technique is Deligne's mixed Hodge structure on the cohomology of any complex algebraic variety. Instead
 of the global analysis on the Riemannian geometry of the complex algebraic variety it will relate to the arithmetic of the variety over finite fields. 

\subsection{Mixed Hodge structure of Deligne}

Motivated by the (then still
unproven) Weil Conjectures and Grothendieck's "yoga of weights", which drew 
cohomological conclusions about complex varieties from the truth of those
conjectures, Deligne in \cite{Deligne-HodgeII, Deligne-HodgeIII} proved the existence of mixed Hodge structures on the cohomology $H^*(M,\Q)$
of a complex algebraic variety $M$. Here we give a quick introduction,
for more details see \cite[\S2.2]{hausel-villegas} and the references therein. 
Deligne's mixed Hodge structure entails two filtrations on the rational cohomology of $M$.  The increasing weight filtration $${0} = W_{-1}  \subseteq W_0   \subseteq \dots \subseteq 
 W_{2j} = H^j(X,\Q)$$
and a decreasing Hodge filtration
$$H^j(X,\C) = F^0  \supseteq  F^1  \supseteq \dots   \supseteq F^m \supseteq    F^{m+1} = {0}.$$ 
We can define mixed Hodge numbers obtained from this two filtrations
by the  following formula:\beq\label{mixedhodgenumbers}h^{p,q;j}(X):=\dim_\C \left(Gr^F_p Gr^W_{p+q}H^j(X)_\C\right).\eeq
 From these numbers we form 
  $$H(M;x,y,t)=\sum_{p,q,k} h^{p,q;k}(M) x^p y^q t^k,$$ the {\em mixed Hodge polynomial}. By virtue of its definition it has the property that the specialization 
$$P(M;t)=H(M;1,1,t)$$ gives  the {\em Poincar\'e polynomial} of $M$.
When $M$ is smooth of dimension $n$ we take another specialization \beq \label{epolynomial} E(M;x,y):=x^ny^nH(1/x,1/y,-1),\eeq the so-called {\em E-polynomial} of a smooth variety $M$. 

Deligne's construction of mixed Hodge structure is complex geometrical: for a smooth variety $M$ 
it is defined by the log geometry of a compactification $\overline{M}$ with normal crossing divisors. In particular a global analytical description, like the Hodge theory of harmonic forms on a smooth complex projective
manifold,  of the mixed Hodge structure on a smooth variety is missing, which causes some difficulty
to find the meaning of mixed Hodge numbers in physical contexts (see the remark after Conjecture~\ref{mirrorderham}).

\subsection{Arithmetic and topological content of the E-polynomial}

The connection of the $E$-polynomial to the arithmetic of the variety is provided by the following theorem of Katz \cite[Appendix]{hausel-villegas}. Here we give an informal version of Katz's result for precise formulation see \cite[Theorem 6.1.2.3, Theorem 2.1.8]{hausel-villegas} :
\begin{theorem}\label{katz} Let $M$ be a smooth quasi-projective variety defined
over $\mathbb Z$ (i.e. given by equations with integer coefficients). Assume that the number of points of $M$ over a finite field $\F_q$, i.e. $$E(q):=\#\{M(\F_q)\}$$ is a polynomial in $q$. Then the $E$-polynomial can be obtained from the count polynomial as follows:
$$E(M;x,y)=E(xy).$$
\end{theorem}
 This theorem is especially useful when we further have $h^{p,q;k}(M)=0$ unless $p+q=k$. In this case we say that the  mixed Hodge structure on $H^*(M)$ is {\em pure}. 
In this case  $$H(M;x,y,t)=(xyt^2)^nE\left(\frac{-1}{xt},\frac{-1}{yt}\right)$$ 
and so the Poincar\'e polynomial can be recovered from the $E$-polynomial as follows $$P(M;t)=H(M;1,1,t)=t^{2n} E\left(\frac{-1}{t},\frac{-1}{t}\right).$$ Examples of varieties with pure MHS on their cohomology include
smooth projective varieties (in this case we get the traditional Hodge structure, which is by definition pure), the moduli space of Higgs bundles $\M_{\Dol}$, the moduli space of flat connections $\M_{\DR}$
on a Riemann surface and Nakajima's quiver varieties.

In general we can define the {\em pure part} of $H(M;x,y,t)$ as $$PH(M;x,y)={\rm Coeff}_{T^0} \left(H(M;xT,yT,tT^{-1})\right).$$ More generally we can define the {\em pure part} of the cohomology of $M$ as $$PH^*(M):=W_n H^n(M)\subset H^*(M),$$ which is a subring $PH^*(M)\subset H^*(M)$  of the cohomology of $M$. For a smooth $M$, the  pure part of $H^*(M)$ 
 is always the image of the cohomology of a smooth 
compactification (see \cite[Corollaire 3.2.17]{Deligne-HodgeII}). 
It is in fact this result which can be used to show that the spaces 
mentioned in the previous paragraph have pure mixed Hodge structure. That is one can prove that they admit a smooth compactification which surjects on cohomology. Prototypes of such compactifications were constructed in   \cite{simpson-hodge} for $\M_{\DR}$ and in \cite{hausel-compact}  for $\M_{\Dol}$.

\section{Applications of mixed Hodge theory}
Using the method sketched in
the previous section the strongest results on cohomology can be achieved when the variety has a pure MHS on its cohomology,  consequently the $E$-polynomial determines the mixed Hodge polynomial, and additionally it is polynomial-count so that Theorem~\ref{katz} gives 
an arithmetic way to determine the $E$-polynomial. This is the
case for Nakajima quiver varieties, where our method gives complete results.

\subsection{Nakajima quiver varieties}

Nakajima quiver varieties are constructed \cite{nakajima-quiver2} by a finite dimensional hyperk\"ahler quotient construction. Here we review
the affine algebraic-geometric version of this construction. 

Let $\Gamma$ be a  quiver (oriented graph) with vertex set $I=\{1,\dots,n\}$ and
edges $E\subset I\times I$.  Let $$\vv=(\vv_1,\dots,\vv_n),\w=(\w_1,\dots,\w_n)\in {\mathbb N}^I$$ be two dimension vectors and
$V_i$ and $W_i$ corresponding complex vector spaces, i.e. $\dim(V_i)=\vv_i$ and $\dim(W_i)=\w_i$. We define the vector spaces  $$\V_{\vv,\w}=\bigoplus_{a\in E} \Hom(V_{t(a)},V_{h(a)})\oplus \bigoplus_{i\in I} \Hom (V_i,W_i)$$ of framed representations of the quiver $\Gamma$ 
, and the  
action $$\rho: \GL(\vv):=\prod_{i\in I} \GL(V_i)\to  \GL(\V_\vv), $$
with derivative $$\varrho: \gl(\vv):=\prod_{i\in I} \gl(V_i)\to  \gl(\V_\vv).$$
The complex moment map  $$\mu: \V_\vv\times \V_\vv^* \to \gl_\vv^* $$ 
of $\rho$ is given at $X\in \gl_\vv$ by  \beq \label{momentmap} 
\langle \mu(v,w),X\rangle=\langle \varrho(X)v,w\rangle. \eeq

 For $\xi=1_\vv \in \gl(\vv)^{\GL(\vv)}$ we define the (always smooth) Nakajima quiver variety by
$${\mathcal M}(\vv,\w)=\mu^{-1}(\xi)/\!/ \GL(\vv)={\rm Spec}\left(\C[\mu^{-1}(\xi)]^{\GL(\vv)}\right)$$ as an affine GIT quotient. Alternatively one can construct the manifold underlying 
${\mathcal M}(\vv,\w)$ as a hyperk\"ahler quotient of $\V_\vv\times \V_\vv^*$ by the maximal compact subgroup $\U(\vv)\subset \GL(\vv)$. This shows that $\M(\vv,\w)$ possesses a hyperk\"ahler metric. The holomorphic symplectic quotient we presented above is the one where 
the arithmetic technique of \S\ref{arithmetic} is applicable. Before we explain that, let us recall the following fundamental theorem  of \cite{nakajima-quiver2} about
the cohomology of these Nakajima quiver varieties:

\begin{theorem}\label{nakajima} Assume that the quiver $\Gamma$ has no edge-loops. Then there is an irreducible representation of the Kac-Moody algebra $\g(\Gamma)$ of highest weight $\w$ on $\oplus_{\vv} H^{mid}({\mathcal M}(\vv,\w))$. In particular the Weyl-Kac
character formula gives the middle Betti numbers of Nakajima quiver
varieties. Furthermore the intersection form on 
$H_c^{mid}({\mathcal M}(\vv,\w))$ is definite, thus $\chi_{L^2}({\mathcal M}(\vv,\w))$ equals the middle Betti number of ${\mathcal M}(\vv,\w)$.
 \end{theorem}

\begin{remark} When $\Gamma$ is an affine Dynkin diagram ${\mathcal M}(\vv,\w)$
could be identified with one of the spaces $M^{k,c_1}_\phi$ of certain Yang-Mills instantons  on a ALE space $X_\Gamma$. 
Kac in \cite{kac-book} explains that 
the Weyl-Kac character formula for an affine Dynkin diagram  has certain modular properties.  This was the line of argument in \cite{vafa-witten} that  \eqref{modular} is a modular form provided Conjecture~\ref{vafawitten} holds.
\end{remark}

In \cite{hausel-betti} a simple Fourier transform technique was found
to enumerate the rational points of $\M(\vv,\w)$ over a finite field 
$\F_q$. The corresponding count function $E(q)$ turned out to be polynomial, and as the mixed Hodge structure is pure on $H^*(\M(\vv,\w))$ the technique of \S\ref{arithmetic} applies in its full strength to give 
a formula for the Betti numbers of the varieties $\M(\vv,\w)$. The result is the following formula from \cite{hausel-betti}:

\begin{theorem}\label{betti-quiver} For any quiver $\Gamma$, 
the Betti numbers of the Nakajima quiver varieties are given by the following generating function, with the notation as in \cite[Theorem 3]{hausel-betti}:
\beq \label{betti} \sum_{{\mathbf v}\in \N^I} P_t(\calM({\mathbf v},{\mathbf w})) t^{-d({\mathbf v},{\mathbf w})}T^{\mathbf v}= \frac{{\displaystyle \sum_{{\mathbf v}\in \N^{I}}} T^{\mathbf v}
 {\displaystyle \sum_{\lambda\in \calP({\mathbf v})}} \frac{
\left( \prod_{(i,j)\in { E}} t^{-2\langle \lambda^i,\lambda^j\rangle}\right)\left(\prod_{i\in {I}} t^{-2\langle \lambda^i,(1^{{\mathbf w}_i})\rangle}\right) }{\prod_{i\in {I}} \left(t^{-2\langle\lambda^i,\lambda^i\rangle}\prod_k \prod_{j=1}^{m_k(\lambda^i)} (1-t^{2j}) \right)}}{{\displaystyle \sum_{{\mathbf v}\in \N^{I}} }T^{\mathbf v}{\displaystyle \sum_{\lambda\in \calP({\mathbf v})}}  \frac{\prod_{(i,j)\in { E}} t^{-2\langle \lambda^i,\lambda^j\rangle}}{\prod_{i\in {I}} \left(t^{-2\langle \lambda^i,\lambda^i\rangle}\prod_k \prod_{j=1}^{m_k(\lambda^i)} (1-t^{2j})\right) }},\eeq\end{theorem}

\begin{remark}  \label{kacconjecture} When $\Gamma$ has no edge-loops Nakajima's Theorem~\ref{nakajima} implies that the right hand side of  \eqref{betti} is a deformation of the 
Weyl-Kac character formula. Simple reasoning gives the same result about the denominator of the right hand side of \eqref{betti} and the Kac denominator. Moreover, the Kac's denominator formula and Hua's formula \cite[Theorem 4.9]{hua}
expressing the denominator of \eqref{betti} as an infinite product implies  a conjecture of Kac, cf. \cite[Corollary 4.10]{hua}.  Namely, if $A_{\Gamma}(\vv,q)$ denotes the number of
absolutely indecomposable representations of $\Gamma$ of dimension vector $\vv$ over the finite field $\F_q$, then it turns out to be a polynomial in $q$ and  Kac's  \cite[Conjecture 1]{kac-quiver} says that  the constant coefficient \beq \label{kac-conj} A_\Gamma(\vv,0)=m_\vv\eeq equals with the multiplicity of the weight $\vv$ in the Kac-Moody algebra $\g(\Gamma)$.  This can be proved, as sketched above and announced in \cite{hausel-betti}, to be  a consequence of \eqref{betti} and the above mentioned results of Nakajima and Hua.
 \end{remark} 

\begin{remark} When the quiver is affine ADE 
and the RHS becomes an infinite product (indications that this can happen is the infinite product in \cite[\S3]{hausel-betti} and the infinite products in the recent \cite{sasaki}) we could get an
alternative proof of the modularity of \eqref{modular} in the Vafa-Witten S-duality conjecture.  \end{remark}

In the remaining part of this survey we will motivate and study another application of the technique in \S\ref{arithmetic}, which will be less powerful as the mixed Hodge structure will fail to be pure, but will also open new interesting directions by the study of this more complicated mixed Hodge structure. 

\subsection{Spaces diffeomorphic to the Hitchin moduli space $\M(C,P_{\U(n)})$} \label{diffeo}

Among the spaces discussed in this paper it is the Hitchin moduli space $\M(C,P_{\U(n)})$ as defined in \S\ref{hitchin} which exhibits perhaps the most plentiful structures many of which are rooted in its hyperk\"ahler quotient origin. In particular there are three distinct complex algebraic variety structures on  $\M(C,P_{\U(n)})$. These
 can be thought  of \cite{simpson-hodge} as the three types of non-Abelian (first) cohomology: Dolbeault, De Rham and Betti, of the Riemann surface $C$. The survey paper \cite{hausel-mln} gives a quick  introduction to these spaces and some of the cohomological
 implications to be discussed  below. 

In this paper the ground field is always $\C$ unless otherwise indicated. Following \cite{hitchin-self, simpson-hodge}  we define a component of the  twisted
$\GL_n=\GL_n(\C)$ Dolbeault cohomology of $C$ as
\bes {\mathcal M}^d_{\Dol}(\GL_n):= 
\left\{\begin{array}{c}\mbox{moduli
space of semistable  rank $n$}\\ \mbox{degree $d$ Hitchin pairs on $C$}\end{array}
\right\}\ees
the $\GL_n$ De Rham cohomology as
\bes {\mathcal M}^d_{\DR}(\GL_n):= 
 \left\{ \begin{array}{c} \mbox{ 
moduli space of flat $\GL_n$-connections} \\\mbox{ on $C\setminus\{ p \}$,
with holonomy $e^{\frac{2\pi i d}{n}}Id$ around $p$}\end{array} \right\} \ees
and the $\GL_n$ Betti cohomology 
\bes{\mathcal M}^d_\B(\GL_n):=\{ A_1,B_1,\dots, A_g,B_g \in \GL_n | \\
 A_1^{-1} B_1^{-1} A_1 B_1 \dots A_g^{-1} B_g^{-1} A_g B_g = e^{\frac{2\pi i d}{n}} Id \}/\!/\GL_n
\ees as a twisted $\GL_n$ character variety of 
$C$. 

When $d=0$ these three varieties are diffeomorphic to the Hitchin moduli space $\M(C,P_{\U(n)})$. However we prefer to consider the twisted versions, when $(d,n)=1$, because then all the varieties are smooth. In this case these three varieties are all diffeomorphic to a twisted version
$\M^d(C,P_{\U(n)})$ of Hitchin moduli space and so to each other. The mixed Hodge structure is pure on $H^*(\M_{\Dol}^d(\GL_n))$ and $H^*(\M_{\DR}^d(\GL_n))$, while it is not pure on $H^*(\M_{B}^d(\GL_n))$. The  mixed Hodge structure are different on $H^*(\M_{\DR}^d(\GL_n))$ and $H^*(\M_{B}^d(\GL_n))$, and so the spaces  $M_{\DR}^d(\GL_n)$ and $\M_{B}^d(\GL_n)$  cannot be isomorphic as complex algebraic varieties. Nevertheless as complex analytic manifolds the Riemann-Hilbert monodromy map \beq \label{RH} \M_{\Dol}^d(\GL_n)\stackrel{RH}{\to} \M_{\DR}^d(\GL_n)\eeq sending a flat connection to its holonomy gives an isomorphism. 

We will also consider the varieties $\M_\Dol^d(\SL_n)$, $\M_\DR^d(\SL_n)$ and $\M_\B^d(\SL_n)$,  which can be defined by replacing $\GL_n$ with $\SL_n$ in the above definitions. Moreover $\M_\Dol^d(\GL_1)$, $\M_\DR^d(\GL_1)$ and $\M_\B^d(\GL_1)$ turn out to be abelian groups. Then $\M_\Dol^d(\GL_1)$, $\M_\DR^d(\GL_1)$ and $\M_\B^d(\GL_1)$, respectively, will act on $\M_\Dol^d(\GL_n)$, $\M_\DR^d(\GL_n)$ and $\M_\B^d(\GL_n)$, respectively, by an appropriate form of tensorization. Finally we denote the corresponding (affine GIT) quotients by $\M_\Dol^d(\PGL_n)$, $\M_\DR^d(\PGL_n)$ and $\M_\B^d(\PGL_n)$. In our case, when $(d,n)=1$, they will turn out to be orbifolds. For more details on the construction of these varieties see 
\cite{hausel-mln}.

In the next section we explain the original motivation to consider the $E$-polynomials of these three complex algebraic varieties.  The motivation is mirror symmetry, and most probably the same $S$-duality we discussed in the Introduction in connection with the Hodge cohomology of the moduli spaces of Yang-Mills instantons in four dimension and magnetic monopoles in three. S-duality ideas relating to mirror symmetry for Hitchin spaces have appeared in the physics literature \cite{bershadsky-etal,kapustin-witten}. 

\subsection{Topological Mirror Test} For our mathematical considerations  the relationship to mirror symmetry stems from the following observation of \cite{hausel-thaddeus-langlands}. It uses the famous {\em Hitchin map} \cite{hitchin-stable}, which makes  the moduli space of Higgs bundles $\M_{\Dol}$ into a completely integrable Hamiltonian system, so that the generic fibers are Abelian varieties.  \begin{theorem}\label{SYZ} In the following diagram  
$$\begin{array}{ccc}
 \M^d_{\Dol}(\PGL_n) &  & \M^d_{\Dol}(\SL_n) \\
\down{\chi_{\PGL_n}} && \down{\chi_{\SL_n}} \\
{\mathcal H}_{\PGL_n} & \cong & {\mathcal H}_{\SL_n}.
\end{array}$$

the generic fibers of the Hitchin maps 
$\chi_{\PGL_n}$ and $\chi_{\SL_n}$ are dual Abelian varieties. 
\end{theorem}

\begin{remark} If we change complex structures and consider
 $\M^d_{\DR}(\PGL_n)$ and 
$\M^d_{\DR}(\SL_n)$, then the Hitchin map on them becomes special Lagrangian fibrations, and consequently the pair of $\M^d_{\DR}(\PGL_n)$ and 
$\M^d_{\DR}(\SL_n)$ satisfies  the requirements of the  SYZ construction \cite{SYZ} for a pair of mirror symmetric
Calabi-Yau manifolds (see \cite{hausel-thaddeus-langlands} and \cite{hausel-thaddeus1} for more details). 
\end{remark}

This motivates the calculation of Hodge numbers of $\M^d_{\DR}(\PGL_n)$ and $\M^d_\DR(\SL_n)$ to see if there is any relationship between them, which one would expect in mirror symmetry. Based on calculations in the $n=2,3$ cases \cite{hausel-thaddeus-langlands} proposed:

\begin{conjecture}\label{mirrorderham}
For all $d,e \in \mathbb Z$, satisfying $(d,n)=(e,n)=1$, 
\bes \Est^{B^e}\left( x,y;\M^d_\DR(\SL_n)\right) 
= \Est^{\hat B^d}\left( x,y;{\mathcal M}^e_\DR(\PGL_n)\right),\ees 
where $B^e$ and $\hat B^d$ are certain gerbes on the corresponding Hitchin spaces and the $E$-polynomials above are stringy $E$-polynomials for orbifolds twisted by the relevant gerbe as defined in \cite{hausel-thaddeus-langlands}.  
\end{conjecture}

Morally, this conjecture should be related to the S-duality considerations of \cite{kapustin-witten} and in turn to the Geometric Langlands Programme of \cite{beilinson-drinfeld}. However the lack of global analytical
interpretation of the  mixed Hodge numbers \eqref{mixedhodgenumbers} 
appearing in Conjecture~\ref{mirrorderham} prevents a straightforward 
physical interpretation. Nevertheless the agreement of certain Hodge numbers for Hitchin spaces for Langlands dual groups is an interesting direction from a purely mathematical point of view. In particular, if we change our focus from $\M_{\DR}$ and $\M_{\Dol}$ to $\M_\B$ we will 
uncover some surprising connections to the representation theory of finite groups of Lie type.

\subsection{Mirror symmetry for finite groups of Lie type}
As $\M_{DR}$ and $\M_\B$ are complex analytically identical via the Riemann-Hilbert map \eqref{RH}, the complex analytical structure of dual special Lagrangian fibrations of Theorem~\ref{SYZ} are present on the pair $\M^d_\B(\SL_n)$ and ${\mathcal M}^e_\B(\PGL_n)$. We might as well try to think of this pair as mirror symmetric in the SYZ picture. The mixed Hodge numbers of $\M_\B$ are however different from the mixed Hodge numbers of $\M_{DR}$ so the corresponding topological mirror test \cite{hausel-mln} will also be different from Conjecture~\ref{mirrorderham}: 

\begin{conjecture}  \label{mirrorbetti} For all $d,e \in \mathbb Z$, satisfying $(d,n)=(e,n)=1$,
\bes\Est^{B^e}\left(x,y,\M^d_\B(\SL_n)\right) 
= \Est^{\hat B^d}\left(x,y,{\mathcal M}^e_\B(\PGL_n)\right). \ees
\end{conjecture}

For this conjecture however there is a powerful arithmetic method to
calculate these $E$-polynomials. Using this technique  we have already managed to check this conjecture \cite{hausel-mln}  when $n$ is a prime and $n=4$.  This arithmetic method is based on the technique explained in \S\ref{arithmetic} and the following character formula from \cite{hausel-villegas}:

\begin{theorem} Let $\G=\SL_n\mbox{ or }\GL_n$, let
 $\G(\F_q)$ be the corresponding finite group of Lie type
\begin{eqnarray*}\begin{array}{l}
E(\sqrt{q},\sqrt{q},\M^d_B(\G))=\#\{{\mathcal M}^d_B(\G({\mathbb F}_q))\}= \sum_{\chi\in Irr(\G({\mathbb F}_q))} 
\frac{|\G({\mathbb F}_q)|^{2g-2}}{\chi (1)^{2g-1}} \chi (\xi^d_n),\end{array} 
\end{eqnarray*}
where the sum is over all irreducible characters of the finite group of 
Lie type $\G(\F_q)$. 
\end{theorem}

This character formula combined with Conjecture~\ref{mirrorbetti}  implies certain relationships between the character tables of $\PGL_n(\F_q)$ and $\SL_n(\F_q)$.  An intriguing way to formulate it is to say that {\em certain differences between the 
character tables
of  $\PGL_n(\F_q)$ and its Langlands dual $\SL_n(\F_q)$ are 
governed by mirror symmetry}. This kind of consideration could be interesting because the character tables of $\PGL_n(\F_q)$ or more generally those of $\GL_n(\F_q)$ have been known for a long time starting with the work of Green \cite{green} in 1955, while the character tables of $\SL_n(\F_q)$ have just recently been completed \cite{bonafe,shoji}. It is especially enjoyable to follow the effect of the mirror symmetry proposal of Conjecture~\ref{mirrorbetti} by comparing the character tables of $\GL_2(\F_q)$ and  $\SL_2(\F_q)$ first calculated a hundred years ago by Jordan \cite{jordan} and Schur \cite{schur}. 

\subsection{Conjectural answer}

Finally, we can put all our observations and conjectures together to state a conjectural answer to the topological side of Problem~\ref{problem}.

As we already noted the mixed Hodge structure on $H^*(\M_\B)$ is not pure. Therefore we are losing information by considering only $E(\M_\B;x,y)$. It turns out that it is interesting to consider the full mixed Hodge polynomial $H(\M_\B;x,y,t)$. When $n=2$ it can be calculated via the explicit description of $H^*(\M_B)$ in 
\cite{hausel-thaddeus-relations}.  We get \cite[Theorem 1.1.3]{hausel-villegas}: \begin{multline*}  H(\M_B(\PGL_2);x,y,t)=\\=
 \frac{(q^2t^3+1)^{2g}}{(q^2t^2-1)(q^2t^4-1)}+
\frac{q^{2g-2}t^{4g-4}(q^2t+1)^{2g}}{(q^2-1)(q^2t^2-1)} -\frac{1}{2}\frac{q^{2g-2}t^{4g-4}(qt+1)^{2g}}{(qt^2-1)(q-1)}-\frac{1}{2}
\frac{ q^{2g-2}t^{4g-4}(qt-1)^{2g}}{(q+1)(qt^2+1)},
\end{multline*} where $q=xy$ and the four terms correspond to the four types of irreducible characters of $\GL(2,\F_q)$.
   When $g=3$ this equals: \begin{multline*}{t}^{12}{q}^{12}+{t}^{12}{q}^{10}+6\,{t}^{11}{q}^{10}+{t}^{12}{q}^{8}+
{t}^{10}{q}^{10} +6\,{t}^{11}{q}^{8}+16\,{t}^{10}{q}^{8}+6\,{t}^{9}{q}^
{8}+{t}^{10}{q}^{6}+{t}^{8}{q}^{8}+26\,{t}^{9}{q}^{6}\\ +16\,{t}^{8}{q}^{
6}+6\,{t}^{7}{q}^{6}+{t}^{8}{q}^{4}+{t}^{6}{q}^{6}+6\,{t}^{7}{q}^{4}+
16\,{t}^{6}{q}^{4}+6\,{t}^{5}{q}^{4}+{t}^{4}{q}^{4}+{t}^{4}{q}^{2}+6\,
{t}^{3}{q}^{2}+{t}^{2}{q}^{2}+1.
\end{multline*}
In particular we see that the pure part is $1+q^2t^4+q^4t^8$. These
terms correspond to the cohomology classes  $1$, $\beta$ and $\beta^2$, and the term $q^6t^{12}$ is not present because by the Newstead relation $\beta^g=\beta^3=0$ holds \cite{hausel-thaddeus-relations}. In particular there is no pure part in the middle $=12$ dimensional cohomology. The same argument holds  for all $g$, which shows that there is no pure part in the middle dimensional cohomology of $\M^1_B(\PGL_2)$. It is however easy to see that the intersection form  on middle cohomology can only be non-trivial on the pure part   and so this implies \cite[Corollary 5.4.1] {hausel-villegas}:
\begin{corollary} The intersection form on $H_{cpt}^*(\M_{B}^1(\PGL_2))$ is trivial. 
\end{corollary}
This gives an alternative proof of \eqref{trivial} as the equation $$\chi_{L^2}(\M^1_B(\SL_2))= \chi_{L^2}(\M^1_B(\PGL_2))$$ is easy to prove. 
Moreover this  approach is more promising to generalize for any $n$. We will offer a conjecture about the 
pure part of the cohomology of $\M^d_\B(\PGL_n)$ below and in turn that
will yield a conjecture for the intersection form on the middle dimensional compactly supported cohomology, answering the topological side of Problem~\ref{problem}. 

To state our conjecture in its full generality we introduce character varieties on Riemann surfaces with $k$ punctures and parabolic type $\mu=(\mu^1,\dots,\mu^k)$ at the punctures, where $\mu^i$ is a partition of $n$.  In other words we fix semisimple conjugacy classes $\calC_1,\dots,\calC_k\subset \GL_n$, which are generic and  have type $\mu$ (in other words $\mu^i_j$ is the multiplicity of the $j$th eigenvalue of a matrix in $\calC_i$). One can prove 
\cite{hausel-aha1}
that there exists generic semisimple conjugacy classes for every type  
$\mu=(\mu^1,\dots,\mu^k)$.    For a generic $\{\calC_1,\dots,\calC_k\}$ of type $\mu$ we define 
\bes \M_\B^\mu:= &\{ A_1,B_1,\dots,A_g,B_g \in \GL_n, C_1\in \calC_1,\dots,C_k\in \calC_k |\\ &  
 [A_1,B_1]  \cdots [A_g,B_g] C_1\cdots C_k=I_n\}/\!/ \GL_n\ees
 as an affine GIT quotient by the diagonal adjoint action of $\GL_n$. 
 The generic choice of the semisimple conjugacy classes implies that $\M^\mu_\B$ is smooth. 
The torus $\GL_1^{2g}$  acts on $\M_\B^\mu$ by multiplying the matrices $A_i$ and $B_i$ by a scalar. We can define the quotient $$\bar{\M}_\B^\mu:=\M_\B^\mu/\!/\GL_1^{2g}$$ as the corresponding $\PGL_n$ character variety. The variety $\bar{\M}_\B^\mu$ is an orbifold. 
 
By studying the Riemann-Hilbert map on the level of cohomologies we are led \cite{hausel-aha1} to consider the crab-shaped quiver  $\Gamma$  associated to  $g$ and $\mu$. Namely, we can put $g$ loops on a central vertex, and $k$ legs of length $l(\mu^j)$. We also equip $\Gamma$ with a dimension vector $\vv$, which has dimension $\sum_{i=1}^l \mu^j_i$ at the $l$th vertex on the $i$th leg. Consider now the number $A_\Gamma(q,\vv)$ of absolutely indecomposable representations of $\Gamma$ of dimension $\vv$ over the finite field $\F_q$. Kac  \cite[Proposition 1.15]{kac-quiver} proved that $A_\Gamma(q,\vv)$ is a polynomial in $q$ with integer coefficients. We have the following conjecture from \cite{hausel-aha1}:

\begin{conjecture}\label{purity}
The pure part of the cohomology of $\bar{\M}_\B^\mu$ is given by 
$$PH(\bar{\M}_\B^\mu,x,y)=(xy)^{d_\mu/2}A_\Gamma(\vv,1/(xy)), $$ where 
$(\Gamma,\vv)$ is the star-shaped quiver and dimension vector given by the  parabolic type $\mu$, and $d_\mu$ is the dimension of ${\M}_\B^\mu$.
\end{conjecture}
This conjecture gives a cohomological interpretation of 
$A_\Gamma(\vv,q)$ and in particular implies that it has non-negative coefficients confirming \cite[Conjecture 2]{kac-quiver} in the case when $\Gamma$ is crab-shaped. When $\mu$ is indivisible Conjecture~\ref{purity} can be proved to follow from the master conjecture in \cite{hausel-letellier-villegas}, which expresses the mixed Hodge polynomials of all the character varieties $\bar{\M}_\B^\mu$ as a generating function generalizing the Cauchy formula for Macdonald polynomials. It also has the following consequence on the topological $L^2$ cohomology $\chi_{L^2}(\bar{\M}^\mu_\B)$ of \eqref{chitopl2}. 

\begin{conjecture}\label{main} The topological $L^2$ cohomology of the manifold $\bar{\M}^\mu_\B$ is given by \beq  \label{g>1} \chi_{L^2}(\bar{\M}_\B^\mu)&=&0, \mbox{ when } g>1\\ \chi_{L^2}(\bar{\M}_\B^\mu)&=&1, \mbox{ when } g=1\label{g=1} \\
 \label{g=0}\chi_{L^2}(\bar{\M}_\B^\mu)&=&m_\vv, \mbox{ when } g=0,\eeq
 where $m_\vv$ is the multiplicity of the weight $\vv$ in the Kac-Moody algebra $\mathfrak{g}(\Gamma)$, which are encoded by the Kac denominator formula for the
 star-shaped quiver $\Gamma$.
\end{conjecture}

When $g>1$ and the parabolic type
is $\mu=((n))$, i.e. we have only one puncture with central conjugacy class, then  one can identify  $\bar{\M}_\B^\mu=\M_\B^d(\PGL_n)$, with some $d$ such that $(d,n)=1$.  In this case \eqref{g>1} 
says that $$\chi_{L^2}\left(\M_\B^d(\PGL_n)\right)=0,$$ which appeared as \cite[Conjecture 4.5.1]{hausel-villegas}. It follows from the mirror symmetry Conjecture~\ref{mirrorderham} that   \bes H_{cpt}^{mid} \left(\M_\B^d(\SL_n)\right)\cong 
H^{mid}_{cpt}\left(\M_\B^d(\PGL_n)\right)\ees and then the intersection forms also agree. This and   \eqref{g=0} then imply that 
\eqref{trivial} holds for any $n$, i.e.  that the intersection form on the
compactly supported cohomology of $\M^d_\B(\SL_n)$ is trivial. This gives a conjectural answer to the topological side of Problem~\ref{problem}. 

When $g=1$ the conjectured \eqref{g=1} follows from Conjecture~\ref{purity} and the observation that the coefficient of $q$ in the $A$-polynomial $A_\Gamma(q)$  for a $g=1$ crab-shaped quiver $\Gamma$  is always $1$. 

When $g=0$ the varieties $\M^\mu_\B=\bar{\M}_\B^\mu$ coincide. Conjecture~\ref{purity} then implies that $$\chi_{L^2}(\M^\mu_\B)=A_\Gamma(\vv,0).$$ Conjecture~\eqref{g=0} is a combination of this and the equality $A_\Gamma(\vv,0)=m_\vv$, that is Kac's \cite[Conjecture 1]{kac-quiver}, which, as discussed in Remark~\ref{kacconjecture},
follows from Theorem~\ref{betti-quiver}.  

Finally one can define $\bar{\M}_{\Dol}^\mu$ the moduli space of stable parabolic $\PGL_n$-Higgs bundles with quasi-parabolic type $\mu^j\in \calP(n)$ and generic weights at the $j$th puncture on the Riemann surface \cite{boden-yokogawa, gothen-etal}. Then one can prove that $\bar{\M}_{\B}^\mu$ is diffeomorphic to $\bar{\M}_{\Dol}^\mu$. Thus Conjecture~\ref{main} also calculates the intersection form on the compactly supported cohomology of the moduli space $\bar{\M}_{\Dol}^\mu$ of stable parabolic $\PGL_n$-Higgs bundles of any rank. 

\begin{example} Consider the genus $0$ Riemann surface $\mathbb{P}^1$ with four punctures. Consider the moduli space  $\M_{toy}$ of stable rank $2$ parabolic Higgs bundles on $\mathbb{P}^1$, with generic parabolic weights on the full parabolic flag at the punctures (see \cite{boden-yokogawa}). This is a complex surface
and the intersection form on $H^2_c(\M_{toy})$ was discussed in \cite[Example 2 for Theorem 7.13]{hausel-compact}. $H^2_c(\M_{toy})$ is $5$ dimensional but $\chi_{L^2}(\M_{toy})$ is only $4$. (The cohomology class of the generic fiber of the Hitchin map is the one in the kernel.)

$\M_{toy}$ is diffeomorphic to the character variety $\bar{\M}^\mu_\B$ where  $g=0$ and
$\mu=((1,1),(1,1),(1,1),(1,1))$. Thus by Conjecture~\ref{main} we should be able to calculate $\chi_{L^2}(\bar{\M}_\B^\mu)$ in terms of the representation theory of the corresponding quiver $\Gamma$. The corresponding quiver $\Gamma$ in this case will be the affine $\tilde{D}_4$ Dynkin diagram, with $\vv=(2,1,1,1,1)$ the minimal positive imaginary root. Its  multiplicity $m_\vv$  in the affine Kac-Moody algebra associated to $\Gamma$ is known to be $4$. Alternatively it is known \cite[Example b to Conjecture 2]{kac-quiver} that $A_\Gamma(\vv,q)=q+4$, which by \eqref{kac-conj} gives $m_\vv=4$. Thus indeed  $\chi_{L^2}(\M_\B^\mu)=m_\vv=4$ checking
\eqref{g=0} in this case  via \cite[Example 2 for Theorem 7.13]{hausel-compact}. 

\end{example}

\end{document}